\documentclass[11pt]{article}
\usepackage{amsmath, amssymb, latexsym}
\usepackage{amsthm} \usepackage{color}
\usepackage{amscd}
\usepackage{amsfonts}

\newcommand{\g}{{\mathfrak g}}
\newcommand{\h}{{\mathfrak h}}
\newcommand{\m}{{\mathfrak m}}

\newcommand{\s}{{\mathfrak s}}
\newcommand{\mt}{\mathfrak t}
\newcommand{\n}{{\mathfrak n}}

\newcommand{\mk}{{\mathfrak k}}
\newcommand{\ma}{{\mathfrak a}}
\newcommand{\mb}{{\mathfrak b}}
\newcommand{\mc}{{\mathfrak c}}

\newcommand{\q}{{\mathfrak q}}
\newcommand{\0}{{\bf 0}}
\newcommand{\C}{{\bf C}}
\newcommand{\R}{{\bf R}}

\newcommand{\Z}{{\bf Z}}

\newcommand{\cL}{{\mathcal{L}}}

\theoremstyle{plain} 
\newtheorem{Th}{\indent\sc Theorem}
\newtheorem{Lm}{\indent\sc Lemma}
\newtheorem{Cr}{\indent\sc Corollary}
\newtheorem{Ps}{\indent\sc Proposition}

\theoremstyle{definition}
\newtheorem{Df}{\indent\sc Definition}
\newtheorem{Rm}{\indent\sc Remark}
\newtheorem{Ex}{\indent\sc Example}
\newtheorem{Nt}{\indent\sc Note}

\newcommand{\im}{\sqrt{-1}\,}

\def\address#1#2{\begingroup
\noindent\parbox[t]{7.8cm}{
\small{\scshape\ignorespaces#1}\par\vskip1ex
\noindent\small{\itshape E-mail address}
\/: #2\par\vskip4ex}\hfill
\endgroup}
\makeatother

\title{
\huge Locally Conformally K\"ahler Structures on Homogeneous Spaces
}

\author{
\textsc{Keizo Hasegawa and Yoshinobu Kamishima}
}
\date{}

\begin{document}

\maketitle


\begin{abstract}
 We will discuss in this paper homogeneous locally
conformally K\"ahler (or shortly homogeneous l.c.K.) manifolds and 
locally homogeneous l.c.K. manifolds from various aspects of study in the field of
l.c.K. geometry. 
 We will provide a survey of known results along with some new results and
observations; in particular we make a complete classification of
$4$-dimensional homogeneous and locally homogeneous l.c.K. manifolds in terms of
Lie algebras. 
\end{abstract}

\section{Introduction} 

A {\em locally conformally K\"ahler structure} 
({\em l.c.K. structure} for short)
on a differentiable manifold $M$
is a Hermitian structure $h$ on $M$ with its associated
fundamental form $\Omega$ satisfying 
$d \Omega = \theta \wedge \Omega$ for some
closed $1$-form $\theta$ (which is so called {\em Lee form}).
Note that a l.c.K. structure $\Omega$ is locally conformally K\"ahler in the sense that
there is a open covering $\{U_i\}$ of $M$ such that $\Omega_i=e^{-\sigma_i} \,\Omega$
is K\"ahler form on $U_i$ for some functions $\sigma_i$, that is, $d \,\Omega_i = 0$.
The condition $d \,\Omega = \theta \wedge \Omega$ is equivalent to the existence of
a global closed $1$-form $\theta$ (Lee form) such that $\theta|U_i = d \sigma_i$.
A differentiable manifold with a l.c.K. structure $(M, h)$  is called a 
{\em locally conformal K\"ahler manifold}.
We know that l.c.K. structure $\Omega$ is globally conformally
K\"ahler (or K\"ahler)
if and only if $\theta$ is exact (or 0 respectively); and a compact l.c.K.
manifold of non-K\"ahler type (i.e. the Lee form is neither $0$ nor exact)
admits no K\"ahler structure compatible with the original complex structure (\cite{V1}).
\smallskip

The inaugural study of l.c.K. manifolds was made by Vaisman in a series of papers
(c.f. \cite{V}, \cite{V1}, \cite{V2}, \cite{V3}) starting with
the cerebrated paper \cite{V}; and since then there have been extensive studies on l.c.K. manifolds 
(c.f. \cite{T}, \cite{DO}, \cite{KO}, \cite{Be}, \cite{GO}).
In this paper we are concerned with l.c.K. structures on
homogeneous and locally homogeneous spaces of Lie groups.
There exist many examples of compact non-K\"ahler
l.c.K. manifolds which are homogeneous or locally homogeneous
spaces of certain Lie groups, such as Hopf surfaces,
Inoue surfaces, Kodaira surfaces, or some class of elliptic
surfaces (c.f. \cite{Be}, \cite{Has}). Their l.c.K. structures are {\em homogeneous}
or {\em locally homogeneous} in the following sense.

\begin{Df}{\rm
{\em A homogeneous locally conformally K\"ahler} ({\em homogeneous l.c.K.} for short)
manifold $M$ is a homogeneous Hermitian manifold with
its homogeneous Hermitian structure $h$, defining a locally
conformally K\"ahler structure $\Omega$ on $M$.
}
\end{Df}

\begin{Df}{\rm
If a simply connected homogeneous l.c.K. manifold $M=G/H$,
where $G$ is a connected Lie group and $H$ a closed
subgroup of $G$, admits a free action of a discrete subgroup
$\Gamma$ of $G$ on the left, then we call a double coset
space $\Gamma \backslash G/H$ 
{\em a locally homogeneous l.c.K.} manifold.
}
\end{Df}


The classification of non-K\"ahler complex surfaces
except those of class VII with $b_2 > 0$ is known:
{\em Kodaira surfaces, 
Inoue surfaces, properly elliptic surfaces of odd type or Hopf surfaces}.
Almost all of these non-K\"ahler complex surfaces, up to small deformations,
admit either homogeneous or locally homogeneous l.c.K. structures
(\cite{T}, \cite{V3}, \cite{GO}, \cite{Be}).
In fact we will see in this paper that except for the class of diagonal Hopf surfaces
with eigenvalues $\lambda_1, \lambda_2 \,(|\lambda_1| \not=|\lambda_2|)$ (which admit
l.c.K. structures of Vaisman type),
we can express, up to small deformations, each of these
complex surfaces as $\Gamma \backslash G/H$,
where $G/H$ is a simply connected $4$-dimensional homogeneous l.c.K. manifold
and $\Gamma$ is a discrete subgroup of $G$ (c.f. \cite{Be}).
\smallskip

It is known (due to Brunella \cite{Br}) that Kato surfaces, 
which are non-K\"ahler complex surfaces of VII with $b_2 > 0$,
also admit l.c.K. structures. And there is a conjecture that
Kato surfaces exhaust all  non-K\"ahler complex surfaces of class VII
with $b_2 >0$. Therefore, it leads to a conjecture that all of non-K\"ahler complex surfaces
admit l.c.K. structures up to small deformations.
\smallskip

Let  $(M, h)$ be a l.c.K. manifold with the Lee form $\theta$. We define the Lee field $\xi$ as
the vector field associated to $\theta$ with respect to the Riemann metric $h$, that is,
$\xi = h^{-1} \theta$. There is an important subclass of l.c.K. manifolds,
called {\em Vaisman manifolds}.

\begin{Df}{\rm
A l.c.K. manifold  $(M, h)$ is of {\em Vaisman type} if
the Lee field $\xi$ is parallel with respect to the Levi-Civita connection of $h$.
}
\end{Df}

For a homogeneous l.c.K. manifold $M=G/H$ with compact isotropy subgroup $H$, 
the Lee field $\xi$ (as an element of the Lie algebra $\g$ of $G$)
is parallel with respect to
the Riemannian connection for $h$ if and only if
$$h(\triangledown_X \xi, Y) = h([X, \xi],Y) - h([\xi,Y],X) + h([Y,X],\xi) =0$$
for all $X,Y \in \g$.
Since the Lee form is closed: $h([Y,X],\xi) =0$, this condition is equivalent to
$$h([\xi, X],Y) + h(X, [\xi,Y])= 0$$
for all $X,Y \in \g$. And this is exactly the case when the Lee field $\xi$ is Killing field. 
It should be also noted that $\xi$ is Killing if and only if $\cL_\xi \Omega=0$ and $\cL_\xi J=0$
for the l.c.K. form $\Omega$ and its compatible complex structure $J$.

\begin{Df}{\rm
We define an exterior differential $d_{\theta}$ on the de Rham complex 
$\Lambda^*(M)$ of a l.c.K. manifold $M$ as 
$$d_\theta: w \rightarrow  - \theta \wedge w + d w,$$
which satisfies $d_\theta^2=0$ for $w \in \Lambda^*(M)$. 
We call $H_\theta^k (M)$ the $k$-th {\em twisted cohomology group}
with respect to $\theta$.}
\end{Df}

It is known  (due to de Le\'on-L\'opez-Marrero-Pard\'on \cite{LLMP}) that 
for a Vaisman manifold $M$, all $H^k_\theta(M)$ vanish.
\smallskip

For a vector space $V$ over $\R$ and a $\g$-module $V$,
we can define $p$-cochains as the $p$-linear
alternating functions on $\g^p$ with values in $V$,
which are $\g$-modules defined by
$$(\gamma f) (x_1, x_2, ...,x_p)= \gamma f(x_1, x_2, ...,x_p) - 
\sum_{i=1}^{p} f(x_1, ...,x_{i-1}, [\gamma, x_i], x_{i+1}, ...,x_p),$$
where $\gamma \in \g$ and $f$ is a $p$-cochain (c.f. \cite{HS}).
The coboundary operator is defined by
\begin{eqnarray*}
(d f) (x_0, x_1, ..., x_p) & = & \sum_{i=0}^{p} (-1)^i x_i f(x_0, ..., \widehat{x_i}, ...,x_p) \\
& + & \sum_{j<k} (-1)^{j+k} f([x_j, x_k], x_0, ..., \widehat{x_j}, ..., \widehat{x_k}, ..., x_p).
\end{eqnarray*}

We are interested in the case where a $\g$-module is defined by
the representation of $\g$ on $\R$, assigning $X \in \g$ to $-\theta(X)$ for
the Lee form $\theta$ on a l.c.K. Lie algebra $\g$.
The corresponding coboundary operator is given by
$$d_\theta: w \rightarrow -\theta \wedge w + d w,$$
and its cohomology group $H_{\theta}^p(\g, \R)$ is nothing but the {\em $p$-th twisted cohomology
group} with respect to the Lee form $\theta$. 
We know (\cite{D}, \cite{HS}) that for a reductive Lie algebra or a nilpotent Lie algebra $\g$, 
all of the cohomology groups $H_{\theta}^p(\g, \R) \; (p \ge 0)$ vanish; and in particular
we have $\Omega= -\theta \wedge \psi + d \psi$ for some $1$-form $\psi$.
Based on this observation we can determine all l.c.K. nilpotent Lie algebras and
l.c.K. reductive Lie algebras (see Section~4). On the other hand, we know (c.f. \cite{R}) that
for a locally homogeneous l.c.K. manifold $\Gamma \backslash G$, where $G$ is
simply connected solvable Lie group with lattice $\Gamma$, there is a canonical injection 
$$H_\theta^k(\g) \hookrightarrow H_\theta^k(\Gamma \backslash G).$$
Hence, for a l.c.K. solvable Lie algebras $(\g, \Omega)$ with
the cohomology class $[\Omega] \not= 0$ in $H_{\theta}^2(\g, \R)$
a compact solvmanifold $\Gamma \backslash G$ is of non-Vaisman type (Note 5).
There exists a l.c.K. Lie algebra  $\g$ of non-Vaisman type of which
all cohomology groups $H_{\theta}^p(\g, \R) \; (p \ge 0)$ vanish (Example 3).
\medskip

There exists a very close relation between l.c.K. structures and Sasaki structures.

\begin{Df} {\rm \hfill
\begin{list}{}{\topsep=0pt \leftmargin=5pt \itemindent=3pt \parsep=0pt \itemsep=3pt}
\item[{$\bullet$}] A {\em contact metric structure} 
$\{\phi, \eta, \widetilde{J}, g\}$ on $M^{2n+1}$
is a contact structure $\phi \, , \phi \wedge (d \phi)^n \not= 0$
with the {\em Reeb field} $\eta \, , i(\eta) \phi = 1, i(\eta) d \phi = 0$,
a $(1, 1)$-tensor $\widetilde{J} \, , \widetilde{J}^2 = -I + \phi \otimes \eta$
and a Riemannian metric $g \, , g(X, Y) = \phi(X) \phi(Y) + d \, \phi (X, \widetilde{J} Y)$.

\item[{$\bullet$}] A {\em Sasaki structure} on $M^{2n+1}$ is
a contact metric structure
$\{\phi, \eta, \widetilde{J}, g\}$ satisfying ${\cal L}_{\eta} g = 0$ (Killing field)
and the integrability of  $J = \widetilde{J}|{\cal D}$ on ${\cal D} = {\rm ker} \, \phi$ (CR-structure).

\item[{$\bullet$}] 
For any Sasaki manifold $N$, its {\em K\"ahler cone} $C(N)$ is defined as
$C(N)=\R_+ \times N$ with the K\"ahler form $\omega=r d r \wedge \phi + \frac{r^2}{2}d \phi$,
where a compatible complex structure $\widehat{J}$ is defined by 
$\widehat{J} \eta = \frac{1}{r} \partial_r$ and $\widehat{J} |{\cal D} = J$.
\end{list}
}
\end{Df}

For any Sasaki manifold $N$ with contact form $\phi$,  we can define
a l.c.K. form $\Omega=\frac{2}{r^2} \omega = \frac{2}{r} d r \wedge \phi + d \phi$;
or taking $t= -2 \,{\rm log}\, r$,
$\Omega= - d t \wedge \phi + d \phi$ on $M=\R \times N$ or $S^1 \times N$,
which is of Vaisman type. We can define a family of complex structures $J$ compatible with
$\Omega$ by
$$J \,\partial_t = b \,\partial_t + (1+b^2) \,\eta, J \eta = - \partial_t - b \, \eta,$$
where $b \in \R$ and the Lee field is $J \eta$.
Conversely, any simply connected complete Vaisman manifold
is of the form $\R \times N$ with l.c.K. structure as above,
where $N$ is a simply connected complete Sasaki manifold.
\medskip

A l.c.K. structure on $M$ may be defined as a K\"ahler structure $\tilde{\omega}$ on 
the universal covering $\tilde{M}$
on which the the fundamental group $\Gamma$ acts homothetically;
that is, for every $\gamma \in \Gamma$, $\gamma^*\tilde{\omega} = \rho(\gamma) \tilde{\omega}$
holds for some positive constant $\rho(\gamma)$. 
\medskip

Let $M=G/H$ be a homogeneous l.c.K. manifold. Then its universal covering
$\tilde{M} =\tilde{G}/\tilde{H}_0$ is also a homogeneous l.c.K. manifold.
Since the Lee form $\tilde{\theta}$ is exact,
$\tilde{\Omega}$ is globally conformal to a K\"ahler structure $\tilde{\omega}$.
The Lie group $\tilde{G}$ acts homothetically
on $\tilde{M}$ from the left, and the fundamental group $\Gamma=\tilde{H}/\tilde{H}_0$
acts on $\tilde{M}$ homothetically from the right. 
Conversely, a K\"ahler structure $\tilde{\omega}$ on $\tilde{M}$
with homothetic action of $\tilde{G}$ from the left and 
$\Gamma$ from the right on $\tilde{M}$ defines a l.c.K. structure on $M$.

\begin{Df}{\rm
Let $M$ be a l.c.K. manifold. Suppose that the universal covering $\tilde{M}$ admits a
{\em K\"ahler potential} $\phi$, which is a real positive function on $\tilde{M}$ such that
$\tilde{\omega}=-\sqrt{-1} \partial \overline{\partial} \phi$
defines a K\"ahler structure on $\tilde{M}$. If the fundamental group $\Gamma$ acts
homothetically on $\phi$, then we call $\phi$ a {\em l.c.K. potential} for $M$.
$\tilde{\omega}$ clearly defines a l.c.K. structure on $M$.}
\end{Df}

We know (due to Ornea-Verbitsky \cite{OV}) that 
a small deformation of a compact l.c.K. manifold with potential
is also a l.c.K. manifold with potential. In other words, l.c.K. structure
with potential is preserved under small deformations. 
\medskip

We know that any primary Hopf surface
can be deformed to a diagonal Hopf surface,
which admits a l.c.K. potential, where a diagonal Hopf surface is 
that of which the contraction map is expressed as a diagonal matrix. Hence we see that
any Hopf surface admits a l.c.K. structure. This argument also holds for
any generalized Hopf manifold of complex dimension $n$
(a compact complex manifold of which the universal
covering is $\C^n -\{0\}$) (c.f. \cite{Has2}); and thus it admits a l.c.K. structure. 

\section{Structure theorems of compact homogeneous l.c.K. manifolds} 
\medskip

We have the following basic result for compact homogeneous l.c.K. manifolds
\cite{HK} (also \cite{ACHK} for an extended form).
\begin{Th} 
A compact homogeneous l.c.K. manifold $M$ is, up to biholomorphism,
isomorphic to a holomorphic principal fiber bundle over a flag
manifold with fiber a $1$-dimensional complex torus $T_\C^1$.
To be more precise, $M$ can be written as a homogeneous space
form $G/H$, where $G$ is a 
compact connected Lie group of holomorphic automorphisms on $M$
which is of the form
$$G = S^1 \times S,$$
where $S$ is a compact simply connected semi-simple Lie group, 
including the connected component $H_0$ of $H$ which is a closed subgroup of $S$. 
$S/H_0$ is a compact simply connected homogeneous Sasaki manifold, 
which is a principal fiber bundle over a flag manifold $S/Q$
with fiber $S^1=Q/H_0$ for some parabolic subgroup $Q$ of $S$
including $H_0$. $M=G/H$ can be expressed as
$$M=S^1 \times_{\Gamma} S/H_0,$$
where $\Gamma = H/H_0$ is a finite abelian group
acting holomorphically on the fiber  $T_{\C}^1$ of the fibration $G/H_0 \rightarrow G/Q$ on the right.
\end{Th}

{\indent\sc Sketch of Proof.} 
Since $G$ is a compact Lie group, it is reductive; and its Lie algebra $\g$ is of the form:
$$\g=\mt + \s,$$
where $\mt$ is the center of $\g$ and $\s$ a semi-simple Lie algebra with $[\g, \g]=\s$.
Since the Lee form $\theta$ is closed but not $0$,
we must have $\theta \in \mt^{*}$. Let $\xi$ be the Lee field with $\theta(\xi)=1$, and
$\eta=J \xi$ (the Reeb field) for the complex structure $J$ with its associated $1$-form
$\phi$ satisfying $\phi(\eta)=1$. We can express $\g$ as
$$\g = <\xi, \eta> +\,\mk,$$
where $<\xi,\eta>$ is the $2$-dimensional subspace of $\g$
generated by $\xi$ and $\eta$ over $\R$, and
$\mk={\rm ker}\,\theta \cap {\rm ker}\,\phi$ with $\mk \supset \h$.
\smallskip

We can see that 
(1) $1 \le \dim \mt \le 2, \, \mt \subset \;<t, \eta>  +\,\h$;
(2) The case ${\rm dim}\,\mt =2$ can be reduced to the case ${\rm dim}\,\mt =1$;
(3) $\xi, \eta$ are infinitesimal automorphisms of $J$ and
infinitesimal isometries (Killing fields) with respect to $h$.
\smallskip

Let ${\mathfrak q}=<\eta>+\h$, then ${\mathfrak q}$ is a Lie subalgebra
of $\s$ and $\h$ is an ideal of ${\mathfrak q}$; in fact we have
${\mathfrak q}=\{X \in \s \,|\, d \phi(X, \s)=0\}$.
Let $S$ and $Q$ be the Lie subgroups of $G$ corresponding to $\s$ and $\q$ respectively.
Then $Q$ is a closed subgroup of $S$ and $H_0$ is a normal subgroup of
$Q$ with $Q/H_0 =S^1$, and $\eta$ generates an $S^1$ action on $S$.
We can see that the Lie subalgebra $<\xi, \eta>$ of $\g$
corresponds to a $2$-dimensional torus $T^2$ of
$G$; $\xi - \sqrt{-1} \eta$ defines a $1$-dimensional
complex torus action on $M=G/H_0$ from the right which is holomorphic
and isometric. We have $M=S^1 \times_{\Gamma} S/H_0$, where $S/H_0 \rightarrow S/Q$
is a principal $S^1$-bundle over the flag manifold $S/Q$ (\cite{BW}, \cite{Bo}); 
and $M \rightarrow S/Q$ is a holomorphic principal fiber
bundle over the flag manifold $S/Q$ with fiber $T_{\C}^1$. 
\bigskip

\begin{Cr}
There exist no compact homogeneous complex l.c.K. manifolds;
in particular, no compact complex paralellizable
manifolds admit their compatible l.c.K. structures.
\end{Cr}

\begin{proof}
Only compact complex Lie groups are complex tori, which can not
act transitively on a compact l.c.K. manifold.
\end{proof}

\section{Homogeneous l.c.K. structures on Lie groups} 

A homogeneous l.c.K. structure on a Lie group $G$ is nothing but
a left invariant l.c.K. structure on $G$. Since $G$ can be expressed
as $\widehat{G}/\Delta$, where \textsf{}$\Delta$ is a finite subgroup of the
center of $\widehat{G}$, $G$ admits a l.c.K. structure $\Omega$
if and only if $\widehat{G}$ admits a l.c.K. structure $\hat{\Omega}$,
or equivalently
the Lie algebra ${\g}$ of $G$ admits a l.c.K. structure $\tilde{\Omega}$
in ${\wedge}\, {\g}^*$.
\medskip

The following result leads to a result of Sawai \cite{S}
on compact locally homogeneous l.c.K. nilmanifolds.
The proof is much simplified but essentially in the same vein
as the original one.

\begin{Th} 
A simply connected (non-abelian) nilpotent Lie group $G$ of dimension $2n$
admits a homogeneous l.c.K. structure
if and only if it is of Heisenberg type: ${\mathcal H}_{2n} =\R \times H_{2n-1}$
with canonical l.c.K. structure, where $H_{2n-1}$ is a $(2n-1)$-dimensional
Heisenberg Lie group.
\end{Th}

\begin{proof}
Let $\g$ be a nilpotent Lie algebra with l.c.K. form $\Omega$.
$\Omega$ is a non-degenerate 2-form such that $d \Omega=\alpha \wedge \Omega$
for some closed 1-form $\alpha$ (Lee form). We know (due to Dixmier \cite{D}) that
there exists a 1-form $\beta$ such that $\Omega= -\alpha \wedge \beta + d \beta$. 

Since $\g$ is a  nilpotent Lie algebra  there exists a non-zero element $B \in [\g,\g] \cap Z(\g)$, 
where $Z(\g)$ is the center of $\g$. Then we have $\Omega(X, B)=- \alpha(X) \beta(B)$ for $X \in \g$.
In particular, since $\Omega$ is non-degenerate we have $\beta(B) \not=0$.
Therefore we may assume, multiplying a positive constant to $B$ if necessary,  that $\beta(B)=1$.
Note that since $\alpha$ is a non-zero closed $1$-form we have $\alpha(B) = 0$. 
Let $A$ be the Lee field  associated to $\alpha$ with respect to $h$, where
$h(X,Y)=\Omega(JX,Y)$ for $X, Y \in \g$. We have $\alpha(Y)=h(A,Y)=\Omega(JA,Y)$
for $Y \in \g$ and $\alpha(A)=1$. Furthermore, taking $\beta + c \alpha$ for $\beta$
with some constant $c \in \R$ if necessary, we can assume that $\alpha, \beta$ satisfy
the condition
$$\alpha(A)=1, \alpha(B)=0, \beta(A)=0, \beta(B)=1.$$

Let $\h$ be the vector subspace of $\g$ generated by $A, B$, and $\n$ the
orthogonal complement of
$\h$ w.r.t. $\Omega$. Since $\Omega$ is non-degenerate, there exist
$X_i, Y_j \in \n$, $i, j=1,...,m$ such that $\n$ is generated by $X_i, Y_j$,
and $d \beta = \sum y_i \wedge x_i$,
where $x_i, y_i$ are the dual forms corresponding to $X_i, Y_i$ 
respectively.
\smallskip

First note that since $\alpha$ is closed we have $A \not\in [\g, \g]$. 
Since $\beta([A,\n]) = - d \beta(A,\n)= -\Omega(A,\n) =0$,
we have $[A,\n] \subset \n$.
We also have $\beta([X_i, Y_j])= - d \beta(X_i, Y_j) = - \Omega (X_i, Y_j) = \delta_{ij}$,
and thus $[X_i, Y_j] = \delta_{ij} B \; {\rm mod}\, \n \;(i, j=1,...,m)$.
Therefore $\n' = <B> + \n$ is an ideal of $\g$. Actually from $\alpha(X)=\Omega(X,B)$
we get $\n'= {\rm ker}\, \alpha$ and $\n =  {\rm ker}\, \alpha \cap {\rm ker}\, \beta$.
\smallskip

Let $\g^{(i)} = [\g, \g^{(i-1)}]$, $\g^{(0)} = \g$. Since $\g$ is nilpotent,
there exists  some positive integer $k$ such that $\g^{(k)} \not= \{0\}$
and $\g^{(k+1)} = \{0\}$.
We will show that $\g^{(k)} = \{B\}$. In fact, any element $Z$ of $\g^{(k)}$
can be written as
$Z = b B + \sum x_i X_i + y_j Y_j \; (b, x_i, y_j \in \R)$.
Then $[Z, Y_i] = x_i B =0 \; {\rm mod} \; \n$,
so $x_i =0$. In the same way we get $y_j=0$;
and thus $\g^{(k)} = \{B\}$.
\smallskip

We will show that $JA=B, JX_i= Y_i, i=1, 2,..., m$.
First note that we have $\alpha(Y)=h(A,Y)=\Omega(JA,Y)$ for any $Y \in \g$.
On the other hand, since $d \beta(B, Y)=0$ we also have $\alpha(Y)=\Omega(B,Y)$
for any $Y \in \g$.
Hence, by non-degeneracy of $\Omega$ we must have $JA=B$. In particular
$J$ preserves $\n$, and thus defines a complex structure on $\n$.
Since we have $h(X_i, X_j)=\Omega(J X_i, X_j)= y_j(JX_i)= \delta^i_j$
and $\Omega(J X_i, A)=\Omega(J X_i, B)=\Omega(J X_i, Y_j)=0$ for $j=1,2,...,m$,
we must have $J X_i=Y_i$ for $i=1, 2,..., m$. 
\smallskip

We can thus consider $\g$ as an extension of $\n'$ by $A$, where $\n'$
is an extension of
$B$ by $\n$:
$$ 0 \rightarrow \n' \rightarrow \g \rightarrow  A \rightarrow 0 $$
$$ 0 \rightarrow B \rightarrow \n' \rightarrow  \n \rightarrow 0 $$
\smallskip

Since $(\n, J)$ is a nilpotent K\"ahler algebra with K\"ahler form $d \beta$, $\n$ must be
abelian (due to Hano \cite{Han}).
Since $B$ belongs to the center of $\g$, $B$ is an infinitesimal automorphism of $J$,
that is, ${\cal L}_B J=0$; and since $A=J B$, $A$ is also
an infinitesimal automorphism of $J$, that is, ${\cal L}_A J=0$.
We see by simple calculation that ${\cal L}_B \Omega = {\cal L}_A \Omega=0$.
Hence, we have  ${\cal L}_B h= {\cal L}_A h=0$, that is, $A, B$ are
Killing fields. In particular, we have $h([A, X], Y) + h(X, [A,Y])=0$ for all $X, Y \in \g$.
We will show that $A$ is in the center of $\g$. Suppose not; then,
since $\n$ is nilpotent, there exists $X_0, Y_0 \in \g$ such that $[A, X_0]=Y_0$ is non-zero
and $[A, Y_0]=0$. Then we have $h([A, X_0], Y_0) + h(X_0, [A,Y_0])= h(Y_0, Y_0)$ is
non-zero, contradicting to the fact that $A$ is Killing field.
\smallskip

We have shown that $\g$ is an extension of  $\h$, which is an abelian ideal
generated by $A, B$,  by the abelian algebra $\n$.

$$ 0 \rightarrow \h \rightarrow \g \rightarrow  \n \rightarrow 0.$$

This completes the proof.
\end{proof}

\begin{Cr} 
A compact nilmanifold $M$ of real dimension $2n$ admits a locally homogeneous
l.c.K. structure if and only if it is of Heisenberg type:
$M = \Gamma \backslash {\mathcal H}_{2n}$, where $\Gamma$ is
a lattice (uniform discrete subgroup) of ${\mathcal H}_{2n}$.
\end{Cr}

\begin{Rm} 
The Lee field on a compact locally homogeneous l.c.K. nilmanifold
$M=\Gamma \backslash {\mathcal H}$ can be irregular for some
lattices $\Gamma$.
\end{Rm}

\begin{Th} 
Let $\g$ be a reductive Lie algebra of dimension $2m$; 
that is, $\g = \mt + \s$,
where $\mt$ is an abelian and $\s$ a semi-simple Lie subalgebra of $\g$
with $\s=[\g,\g]$.
Then $\g$ admits a l.c.K. structure if and only if ${\rm dim}\, \mt = 1$
and ${\rm rank} \, \s =1$. In particular a compact Lie group admits
a homogeneous l.c.K. structure if and only if it is $U(2)$,
$S^1 \times SU(2) \cong S^1 \times Sp(1)$, or $S^1 \times SO(3)$;
and any homogeneous l.c.K. structure on a compact Lie group is of Vaisman type.
\end{Th}

\begin{proof}
Suppose that $\g$ admits a l.c.K. structure $\Omega$.
Since we have  $\h=\{0\}$, $\eta \in \s$ and thus ${\rm dim}\, \mt =1$.
If we apply the proof of Theorem 1 for the case $\h=\{0\}$,
we see that $\q =<\eta>=\{V \in \s |\, [\eta, V]=0\}$;
and thus ${\rm rank} \, \s =1$ (c.f. \cite{BW}). We know all of the reductive Lie algebras
$\g=\mt + \s$ with ${\rm dim}\, \mt = 1$ and ${\rm rank} \, \s =1$:
$\R \oplus {\mathfrak sl}(2, \R)$ and 
${\mathfrak u}(2)=\R \oplus {\mathfrak su}(2)=\R \oplus {\mathfrak so}(3)$.
We show that all homogeneous l.c.K. structures on ${\mathfrak u}(2)$ are the ones 
we obtained in Proposition 3: $\Omega=-\theta \wedge \phi+ d \phi$; and
they are all of Vaisman type. In fact, any l.c.K. form $\Omega'$ is of the form
$$\Omega'=-\theta \wedge \psi + d \psi,$$
where we can set $\theta=t$ and $\psi=a x+b y+c z \,(a, b, c \in \R)$; and thus
$d \psi=-(a \,y \wedge z+ b \,z \wedge x + c \,x \wedge y)$.
For the complex structure $J_{\delta}$ in Proposition~2, we denote by $A$ the
$4\times4$-matrix determined by $h'(U,V)=\Omega'(J_{\delta} U,V)$
for $U,V=T,X,Y,Z$. By the condition that $A$ is a positive-definite symmetric matrix,
we can see by calculation that $b=c=0$; and thus $A=a I_4$. Hence
$\Omega'$ is equal to the original $\Omega$ up to constant multiplication.
\end{proof}

\section{Compact homogeneous and locally homogeneous 
l.c.K. manifolds of complex dimension $\bf 2$} 
\label{sec:5}

We know (c.f. \cite{V1}, \cite{GO}, \cite{Be}) that there is a class of Hopf surfaces
which admit homogeneous l.c.K. structures. We can show,
applying Theorem 1, that the only compact homogeneous l.c.K.
manifolds of complex dimension $2$ are Hopf surfaces of homogeneous
type (Theorem~2). We first determine, recalling a result of Sasaki \cite{SS},
all homogeneous complex structures on $G=S^1 \times SU(2)$,
or equivalently all complex structures on the Lie algebra
$\g={\mathfrak u}(2)$.
\medskip

\begin{Ps} 
Let $\g={\mathfrak u}(2)=\R \oplus \s {\mathfrak u}(2)$ be a reductive
Lie algebra with basis $\{T, X,Y,Z\}$ of $\g$, where $T$ is a generator
of the center $\R$ of $\g$, and
$$
X=\frac{1}{2}\left(
\begin{array}{cc}
\im & 0\\
0 & -\im
\end{array}
\right),
\;
Y=\frac{1}{2}\left(
\begin{array}{cc}
0 & \im\\
\im & 0
\end{array}
\right), 
\;
Z=\frac{1}{2}\left(
\begin{array}{cc}
0 & -1\\
1 & 0
\end{array}
\right)
$$
such that non-vanishing bracket multiplications are given by
$$[X,Y]=Z,\;[Y,Z]=X,\;[Z,X]=Y.$$
Then $\g$ admits a family of complex structures $J_\delta, \delta=c+\im d$
defined by
$$J_\delta(T-d X)=c X,\; J_\delta(c X)=-(T-d X),\; J_\delta Y= \pm Z,
\; J_\delta Z = \mp Y.$$
Conversely, the above family of complex structures exhaust all homogeneous
complex structures on $\g$.
\end{Ps}

\begin{proof} 
Let $\g_{\C}=\g {\mathfrak l} (2,\C)=\C+\s {\mathfrak l} (2,\C)$
be the complexficaion of $\g$,
which has a basis ${\mathfrak b}_{\C} =\{T,U,V,W\}$, where
$$U=\frac{1}{2}\left(
\begin{array}{cc}
-1 & 0\\
0 & 1
\end{array}
\right),
\;
V=\frac{1}{2}\left(
\begin{array}{cc}
0 & 0\\
1 & 0
\end{array}
\right),
\;
W=\frac{1}{2}\left(
\begin{array}{cc}
0 & 1\\
0 & 0
\end{array}
\right)
$$
with the bracket multiplication defined by
$$[U,V]=V,\;[U,W]=-W,\;[V,W]= \frac{1}{2} U.$$
Here we have
$$U=\im X,\;
V=\frac{1}{2}(Z-\im Y),\;
W=-\frac{1}{2}(Z+\im Y),$$
and their conjugations given by
$$\overline{T}=T,\; \overline{U}=-U,\; \overline{V}=-W,\;
\overline{W}=-V.
$$

We know that there is a one to one correspondence between complex
structures $J$ and complex subalgebras $\h$ such that
$\g_{\C} =\h + \overline{\h}$
and $\h \cap \overline{\h}=\{0\}$. 
Let $\ma$ be the subalgebra
of $\g_{\C}$ generated by $T$ and $\mb$ the subalgebra of
$\g_{\C}$ generated by $U,V,W$, then we have
$$\g_{\C} = \ma \oplus \mb$$
where $\ma=<T>_{\C}, \mb=<U, V, W>_{\C}$.
Let $\pi$ be the projection $\pi: \g_{\C} \rightarrow \mb$ and
$\mc$ the image of $\h$ by $\pi$, then we have 
$$\mb=\mc+\overline{\mc},$$
and ${\rm dim}\, \mc \cap \overline{\mc} =1$.
We can set a basis $\eta$ of $\h$ as
$\eta=\{P+Q, R\}\; (P \in \ma, Q,R \in \mb)$
such that $Q \in \mc \cap \overline{\mc}$ and
$\gamma =\{Q,R\}$ is a basis of $\mc$:
$$\h=<P+Q, R>_{\C},\; \mc=<Q, R>_{\C}.$$
Furthermore, we can assume that $Q+\overline{Q}=0$
so that $Q$ is of the form $aU+bV+\overline{b}W \,(a \in \R, b \in \C)$.
\smallskip

We first consider the case where $R=qV+rW \, (q,r \in \C)$. 
Since we have $[\g_{\C}, \g_{\C}]=\mb$, there is some $\alpha \in \C$
such that $[Q,R]=\alpha R$.
We see by simple calculation that if $b \not=0$,
then $q=sb, r=s \overline{b}$
for some non zero constant $s \in \C$. But then $\overline{R}=-\frac{\bar{s}}{s} R$,
contradicting to the fact that $\beta=\{Q, R,\overline{R}\}$ consists
a basis of $\mb$:
$$\mb=<Q, R, \overline{R}>_{\C}$$
Hence we have $b=0$, and $q \not=0, r =0$ with
$\alpha=a$ or $q =0, r \not=0$ with $\alpha=-a$.
Therefore we can take, as a basis of $\h$, $\eta=\{T+\delta U, V\}$ or
$\{T+\delta U, W\}$ with $\delta=c+\im d \in \C$:
$$\h=<T+\delta U, V>_{\C} \;{\rm or}\; <T+\delta U, W>_{\C}.$$
It should be noted that
the latter defines a conjugate complex structure of the former,
which are not equivalent but define biholomorphic
complex structures on its associated Lie group $G$.
\smallskip

In the case where $R=pU+qV+rW,\,p, q, r \in \C$ with $p \not=0$,
we show that there exists an automorphism $\widehat{\phi}$ on $\g_{\C}$ 
which maps $\h_0$ to $\h$, preserving the conjugation, where
$\h_{0}$ is a subalgebra of $\g_{\C}$ of the first type with $p=0$.
As in the first case, we must have $[Q,R]=\eta R$ for some non zero
constant $\eta \in \C$. We may assume that $p=1$.
We see, by simple calculation that
$b, q, r  \not= 0$ and
$$(a-\eta)q=b, (a+\eta)r=\overline{b},$$
from which we get
$$a^2+|b|^2=\eta^2 \;(\eta \in \R),$$
and
$$|q|^2-|r|^2 = \frac{4 a \eta}{|b|^2}.$$
Then an automorphism $\phi$ on $\mb$ defined by 
$$\phi(U)=\frac{1}{\eta}Q,\, \phi(V)=\frac{|b|}{2 \eta} R, \, \phi(W)=-\frac{|b|}{2 \eta} \overline{R},$$
extends to the automorphism $\widehat{\phi}$
on $\g_{\C}$ which satisfies the required condition.
\end{proof}
\medskip

\begin{Ps} 
Let $G=S^1 \times SU(2)$ (which is, as is well known,
diffeomorphic to $S^1 \times S^3$). Then
all homogeneous complex structures on
$G$ admit their compatible homogeneous l.c.K. structures,
defining a primary Hopf surfaces $S_{\lambda}$ which are compact
quotient spaces of the form
$W /\Gamma_{\lambda}$, where $W=\C^2 \backslash \{\0\}$
and $\Gamma_{\lambda}$
is a cyclic group of holomorphic automorphisms on $W$ generated by
a contraction $f: (z_1,z_2) \rightarrow
(\lambda z_1, \lambda z_2)$ with $ |\lambda| \not=0, 1$. 
Furthermore, all of these l.c.K. structures are of Vaisman type. 
A secondary Hopf surface with homogeneous l.c.K. structure is obtained as
a finite quotient of a primary Hopf surface: $S^1\times_{\Z_m} SU(2)$.
\end{Ps}

\begin{proof} 
We consider the following canonical diffeomorphism $\Phi_\delta$,
which turns out to be biholomorphic for each homogeneous complex
structure $J_\delta$ on $\g$
and $\lambda_\delta$:
$$\Phi_\delta: \R \times SU(2) \longrightarrow W$$
defined by
$$(t,z_1,z_2) \longrightarrow (\lambda_\delta^t z_1, \lambda_\delta^t z_2),$$
where $SU(2)$ is identified with
$S^3=\{(z_1,z_2) \in \C^2 \,|\;|z_1|^2+|z_2|^2=1\}$
by the correspondence:
$$
\left(
\begin{array}{cc}
z_1 & -\overline{z}_2\\
z_2 & \overline{z}_1
\end{array}
\right)
\longleftrightarrow (z_1,z_2),
$$
and $\lambda_\delta=e^{c+\im d}$.
Then we see that $\Phi_\delta$ is a biholomorphic map.
It is now clear that $\Phi_\delta$ induces a biholomorphism
between $G=S^1 \times SU(2)$ with homogeneous complex
structure $J_\delta$ and a primary Hopf surface 
$S_{\lambda_\delta}=W/\Gamma_{\lambda_\delta}$.
\smallskip

Let $t, x, y, z \in \g^*$ be the Maurer-Cartan forms corresponding to
$T, X,Y,Z \in \g$  in Proposition 2. Then we have
$$d z=- x \wedge y,\, d x=-y \wedge z,\, d y=- z \wedge x,$$
and
$$\Omega=- \theta \wedge \phi + d \phi,$$
where $\theta=t,\, \phi=\frac{1}{c} x,$
defines a l.c.K. form on $\g$ for the complex structure $J_{\delta}$ in 
Proposition 2.
Note that we have the Lee field $\xi=T-\frac{d}{c} \eta$, which
is irregular for an irrational $\frac{d}{c}$ while
the Reeb field $\eta=c X$, which is always regular. The Lee field $\xi$ is
a Killing field, since we have
$$h([\xi, U],V)+h(U,[\xi,V])= -d (h([X, U],V)+h(U,[X,V]))=0$$
for all $U,V \in \g$. Hence $(G; \Omega, J_\delta)$ is of Vaisman type.
\smallskip

A secondary Hopf surface with homogeneous l.c.K. structure can be obtained
as a quotient space of a primary Hopf surface $S_{\lambda_\delta}$ by
some finite subgroup of $G$.
For instance, $U(2)$ is a quotient Lie group of $G$ by
the central subgroup $\Z_2=\{(1, I), (-1, -I)\}$. In general
we have a secondary Hopf surface 
$G/\Z_m = S^1 \times _{\Z_m} SU(2),$
where ${\bf Z}_m$ is a
finite cyclic subgroup of $G$ generated by $c$:
\[
c= (\xi, \tau),\quad \tau = \left(
\begin{array}{@{}cc@{}}
\xi^{-1} & 0\\
0 & \xi
\end{array}
\right)\!,\enspace \xi^m = 1,
\]
with homogeneous l.c.K. structures induced from those on $G$ by the averaging
method (c.f. \cite{Has}). A (primary or secondary) Hopf surface defined as above
is called a {\em Hopf surface of homogeneous type}, which is a holomorphic principal
bundle over a $1$-dimensional projective space $\C P^1$ with fiber a $1$-dimensional
complex torus $T_{\C}^1$.
\end{proof}
\medskip

\begin{Th} 
Only compact homogeneous l.c.K. manifolds of complex dimension $2$
are Hopf surfaces of homogeneous type (up to biholomorphism).
\end{Th}

\begin{proof} 
It is sufficient to show that any compact homogeneous l.c.K.
manifold $M$ of complex dimension $2$ is a Hopf surface of homogeneous type as defined
in Proposition 2. As we have seen in Theorem~1, a compact homogeneous
l.c.K. manifold $M$ of complex dimension $2$ can be expressed as
$S^1 \times_{\Gamma} S$, where $S$ is a compact homogeneous contact
manifold of real dimension $3$ which admits a Hopf fibration over
$\C P^1$ with fiber $S^1$, and $\Gamma$ is a finite abelian group acting on the fiber $T_{\C}^1$
of the fibration $M \rightarrow \C P^1$.
These are exactly Hopf surfaces with homogeneous
l.c.K. structures as defined in Proposition 2.
Conversely a Hopf surface of homogeneous type admits a homogeneous l.c.K. structure as
defined in Proposition~2. 
\end{proof}
\bigskip

We classify, up to biholomorphism, all compact homogeneous and locally homogeneous 
l.c.K. manifolds of complex dimension $2$. To be more precise,
there exist three types of compact locally homogeneous complex manifolds:
(1) homogeneous complex manifold $G/H$, (2) locally homogeneous complex manifold
$\Gamma \backslash G$, where $G$
is a simply connected homogeneous complex unimodular Lie group
with a lattice $\Gamma$, (3) locally homogeneous
complex manifold $\Gamma \backslash G/H$, where $G/H$ is a simply
connected homogeneous complex manifold
with a non-trivial closed subgroup $H$ and a discrete subgroup $\Gamma$ of $G$.
\smallskip

We first classify all unimodular
Lie groups of real dimension $4$; then pick up those admitting lattices.

\begin{Ps} 
Let $\mathfrak n$ be the real nilpotent Lie algebra of dimension $3$,
which has a basis $\{X, Y, Z\}$ with bracket multiplication defined by
$[X,Y]=[X,Z]=0$ and $[Y, Z] = -X$. Then, there are two classes
of unimodular solvable Lie algebras $\g$ of
the form $\n \rtimes \R$ (semi-direct sum) with $\n =[\g, \g]$,
where the adjoint representation is given by the following:
\begin{list}{}{\leftmargin=15pt \itemindent=-5pt \parsep=7pt}
\item{(1)} $ad_W X=0, ad_W Y = -Z, ad_W Z = Y$,
\item{(2)} $ad_W X=0, ad_W Y = Y, ad_W Z = -Z$,
\end{list}
where $W$ is a generator of $\R$ and $ad_W V=[W,V]$.
\end{Ps}

\begin{proof} 
Let $ad_W Y$ and $ad_W Z$ is given by
$$ad_W Y = a Y+b Z+p X,\;  ad_W Z = c Y+d Z+q X.$$

By Jacobi identity, we have
$$[ad_W Y,Z]+[Y,ad_W Z]=ad_W [Y,Z].$$

It follows that $ad_W X = (a+d) X$. Since $\g$ is unimodular, $a+d=0$;
and in particular, we have $ad_W X =0$, that is, $X$ generates
the center of $\g$.
Therefore, $\g$ is given by the following bracket multiplication:
$$[Y,Z]=-X, [W,X]=0, [W,Y]=a Y+b Z, [W,Z]= c Y+(-a) Z.$$

In case  $a=0$, since $[\g,\g]=\n$, we must have $b \not=0, c \not=0$.
If we put $W' = \frac{1}{c} W, X'=-\frac{b}{c} X,Y'=Y, Z'=-\frac{b}{c} Z$,
then we get
$$[Y',Z']=-X', [W',Y']=-Z', [W',Z']=Y'.$$

In case $a \not=0$, the coefficient matrix  has the determinant
$-a^2-bc \not= 0$ and trace $0$.  If $a^2+bc> 0$, then it is
diagonalizable with eigenvalues $\alpha, \beta \, (\alpha+\beta=0)$.
Hence we can assume $\g$ has the form
$$[Y,Z]=-X, [W,Y]=Y, [W, Z]=-Z.$$

If $a^2+bc < 0$, then we can assume that  $\g$ has the form
 $$[Y,Z]=-X, [W,Y]=Y+b Z, [W,Z]=c Y-Z,$$
where $b, c \not= 0$ and $1+bc < 0$.

If we put $Y'=Y-\frac{1}{c} Z$, then we have
$$[Y' , Z]=-X, \, [W, Y']=\frac{bc+1}{c} Z, \, [W, Z]=c Y'$$

If we put
$Z' = - \frac{bc+1}{c} Z, \,W'= \frac{1}{c} W, \, X'= - \frac{bc+1}{c} X$,
then we get
$$[Y',Z']=-X', \, [W',Y']=-Z', \, [W',Z']=Y'$$
\end{proof}

\begin{Nt} 
{\rm A simply connected solvable Lie group associated  with
the unimodular solvable Lie algebras $\g$ of type (1), or of type (2)
admits a lattice, defining a secondary Kodaira surface, or
an Inoue surface of type $S^+$ respectively. In particular,
any simply connected solvable Lie group of dimension $4$
corresponding to a unimodular solvable Lie algebra $\g$ of
the form $\n \rtimes \R$  with $\n =[\g, \g]$ admits a lattice.
}
\end{Nt}

\begin{Ps} 
There are six classes of unimodular solvable Lie algebras $\g$ of the form
$\R^3 \rtimes \R$, where the adjoint representation is given by the following:
$$ad_W X_i=\sum_{j=1}^{3} a_{ij} X_j \; ,i=1,2,3,$$
where $\{X_1,X_2,X_3\}$ is a basis of $\R^3$, $W$ a generator of $\R$,
and $A=(a_{ij})$ a $3 \times 3$ real matrix with ${\rm Tr} A=0$.
Taking a suitable basis and a generator, we can classify $A$ into six types
according to its eigenvalues:

\begin{list}{}{\leftmargin=15pt \itemindent=-5pt \parsep=7pt}
\item{(3)} all the eigenvalues are zero

$$
(i)  \left(
\begin{array}{ccc}
0 & 1 & 0\\
0 & 0 & 1\\
0 & 0 & 0
\end{array}
\right), \hspace{7mm}
(ii)  \left(
\begin{array}{ccc}
0 & 1 & 0\\
0 & 0 & 0\\
0 & 0 & 0
\end{array}
\right).
$$

\item{(4)} only one of the eigenvalues is zero
$$ \left(
\begin{array}{ccc}
0 & 0 & 0\\
0 & a & 0\\
0 & 0 & -a
\end{array}
\right), \;  a \in \R \, (a \not=0).
$$

\item{(5)} zero and  pure imaginary complex eigenvalues 

$$ \left(
\begin{array}{ccc}
0 & 0 & 0\\
0 & 0 & -b\\
0 & b & 0
\end{array}
\right), \;  b \in \R \, (b \not=0).
$$

\item{(6)} three distinct real eigenvalues

$$ \left(
\begin{array}{ccc}
a & 0 & 0\\
0 & b & 0\\
0 & 0 & -(a+b)
\end{array}
\right), \; a, b \in \R \, (a,b \not=0).
$$

\item{(7)} a single and double eigenvalues
$$
(i)  \left(
\begin{array}{ccc}
-2a & 0 & 0\\
0 & a & 0\\
0 & 0 & a
\end{array}
\right), \hspace{7mm}
(ii)  \left(
\begin{array}{ccc}
-2a & 0 & 0\\ 
0 & a & 1\\
0 & 0 & a
\end{array}
\right), \; a \in \R \, (a  \not=0).
$$

\item{(8)} one real and non-real complex eigenvalues

$$ \left(
\begin{array}{ccc}
-2 a & 0 & 0\\
0 & a &- b\\
0 & b & a
\end{array}
\right), \; a, b \in \R \, (a,b \not=0).
$$

\end{list}

\end{Ps}

\begin{Nt} 
{\rm All simply connected solvable Lie group associated  with
the unimodular solvable Lie algebras $\g$ of the above types
except (7) admits  lattices, defining compact solvmanifolds of
dimension $4$. The unimodular solvable Lie algebra $\g$ of
type (3)(ii), (5) and (8) with suitable $a, b$ is corresponding
to Kodaira surface, Hyperelliptic surface and Inoue surface of type
$S^0$ respectively.

We can see, from the following lemma that a simply connected
solvable Lie group associated  with the unimodular
solvable Lie algebras $\g$ of type (7) does not admit any lattice.
}
\end{Nt}

\begin{Lm} 
Let $\Phi(t)$ be a polynomial of the form
$\Phi(t)=t^3-m t^2+n t-1 \; (m, n \in \Z)$.
Then, it has a real double root $a$ if and only if $a=1$ or $-1$
for which $\Phi(t)=t^3-3 t^2+3 t-1$ or $\Phi(t)=t^3+t^2-t-1$ respectively.
\end{Lm}

\begin{proof} 
Assume that $\Phi(t)$ has a double root $a$ and another root $b$.
Then we have that
$a^2 b=1, 2 a+b=m, a^2+2 a b=n$, from which we deduce that 
$m a^2-2 n a+3=0, 3 a^2-2 m a+n=0$; and thus $2 (m^2-3 n) a=m n-9$.
If $m^2=3 n$, then $m=n=3$ and $a=1$. If $m^2 \not= 3 n$, then we have that
$a=\frac{m n-9}{2 (m^2-3 n)}$, which is a rational number. Since we have
that $2 a+\frac{1}{a^2}=m \in \Z$, $a$ must be $1$ or $-1$.
\smallskip

A lattice $\Gamma$ of a simply connected solvable Lie group associated
with the unimodular solvable Lie algebras $\g$ of type
$\R^3 \rtimes \R$ is of the form  $\Z^3 \rtimes \Z$, where the action
$ \phi: \Z \rightarrow {\rm Aut}(\Z^3)$ is determined by
$\phi(1) = A \in {\rm SL}(3, \Z)$;
and the characteristic polynomial $\Phi(t)$ of $A$ is of the form
$\Phi(t)=t^3-m t^2+n t-1 \; (m, n \in \Z)$. According to the above lemma,
$\Phi(t)$ can have a double root $a$ if and only if $a=1$ or $-1$.
\end{proof}

\begin{Th} 
There are ten classes of unimodular Lie algebras of dimension $4$;
eight classes of solvable Lie algebras obtained in the proposition 3 and 4, 
and two classes of reductive Lie algebras: 
${\mathfrak gl}(2, \R) = \R \oplus {\mathfrak sl}(2, \R)$
and ${\mathfrak u}(2) = \R \oplus {\mathfrak su}(2)$.
Their associated simply connected reductive Lie groups admit lattices,
defining a properly elliptic surface and a Hopf surface respectively.
\end{Th}

\begin{proof} 
Applying Levi decomposition, a Lie algebra of dimension $4$ is
either solvable or reductive of the form $\R \oplus {\mathfrak s}$,
where $\mathfrak s$ is a simple Lie algebra, which is either
$ {\mathfrak sl}(2, \R)$ or ${\mathfrak su}(2)$.
\end{proof}
\medskip

We will see that most of non-K\"ahler complex surfaces of the form
$\Gamma\backslash G$ with a unimodular Lie group $G$ having a lattice
$\Gamma$ admit locally homogeneous l.c.K. structures (\cite{Has}).
In the following list the Lie algebra $\g$ is generated by $X,Y,Z,W$
with the specified brackets multiplication.
\smallskip

\begin{list}{}{\leftmargin=15pt \itemindent=-5pt \parsep=7pt}
\item{(1)} Primary Kodaira surface:
$[X,Y]=-Z$, and all other brackets vanish.

\item{(2)} Secondary Kodaira surface:
$[X,Y]=-Z,\, [W,X]=-Y,\, [W,Y]=X$, and all other brackets vanish.

\item{(3)} Inoue surface $S^{\pm}$:
$[Y,Z]=-X,\, [W,Y]=Y,\, [W,Z]=-Z$, and all other brackets vanish.

\item{(4)} Inoue surface $S^0$:
$[W,X]=-\frac{1}{2} X-b Y,\, [W,Y]=b X -\frac{1}{2} Y,\, [W,Z]= Z$,
and all other brackets vanish.

\item{(5)} Properly elliptic surface:
$[X,Y]=-Z,\, [Z,X]=Y,\, [Z,Y]=-X$, and all other brackets vanish.

\item{(6)} Hopf surface:
$[X,Y]=-Z,\, [Z,X]=-Y,\, [Z,Y]=X$,  and all other brackets vanish.
\end{list}

For all cases, we have a homogeneous complex
structure defined by
$$J Y=X, J X=-Y, J W=Z, J Z=-W,$$
and its compatible l.c.K. form
$\Omega =x \wedge y + z \wedge w$ with
the Lee form $\theta=w$, where $x, y, z, w$ are the Maurer-Cartan
forms corresponding to $X,Y,Z,W$ respectively. 

\begin{Nt} 
For Inoue surfaces of type $S^+$, we have other homogeneous
complex structures:
$$J Y=X, J X=-Y, J W=Z+q X, J Z=-W+q Y,$$
with non-zero real number $q$,
for which there exist no compatible l.c.K. structures (due to Belgun \cite{Be}).
\end{Nt}

\begin{Nt} 
For properly elliptic surface and Hopf surfaces, 
we have other homogeneous complex structures
$$J Y=X, J X=-Y, J W=b W+(b^2+1) Z, J Z=-W- b Z,$$
with no-zero real number $b$, for all of which $\Omega$ defines a compatible l.c.K. structure.
Any locally homogeneous l.c.K. Hopf surface is of Vaisman type, while some of
locally homogeneous l.c.K. properly elliptic surfaces are of non-Vaisman type (Example 3).
\end{Nt}

\begin{Nt} 
\hfill
\begin{list}{}{\topsep=0pt \leftmargin=12pt \itemindent=0pt \parsep=0pt \itemsep=3pt}
\item[{$\bullet$}] 
For secondary Kodaira surface,  the Lee filed $\xi = W$, and the bracket multiplication
is given by $[X,Y]=-Z,\, [W,X]=-Y,\, [W,Y]=X$.
We get by simple calculation,
$$h(\triangledown_U W, V) = h([W, U],Y) + h(U, [W, V]) = 0$$
for any $U, V \in \g$. Hence $\xi$ is a Killing field.
It is also easy to check $\Omega=- w \wedge z + d z$. 

\item[{$\bullet$}] 
For Inoue surface $S^{\pm}$,  the Lee filed $\xi = W$, and the bracket multiplication
is given by $[Y,Z]=-X,\, [W,Y]=Y,\, [W,Z]=-Z$. We get by simple calculation,
$$h(\triangledown_Z W, Z) = h([W, Z],Z) + h(Z, [W, Z]) = -2 h(Z, Z) \not= 0.$$
Hence $\xi$ is not a Killing field.
It is also easy to check that there is no invariant $1$-form $v$ such that
$\Omega=- w \wedge v + d v$.
\end{list}

\end{Nt}

\begin{Th} 
Only compact locally homogeneous l.c.K. manifolds
$\Gamma\backslash G/H$ of real dimension $4$
with non-trivial closed subgroup $H$ are Inoue surfaces of type $S^-$ (Example 1),
Hopf surfaces of locally homogeneous type (Example 2) and properly elliptic surfaces.
\end{Th}

\begin{proof} 
Let $M=G/H$ be a simply connected homogeneous l.c.K. manifold of
dimension $4$,
where $G$ is a connected Lie group with
a non-trivial compact subgroup $H$. We have the Levi-decomposition
$G=R \cdot S$, where $R$ is the radical of $G$ which can be written as
$R=F \cdot T \, (F \cap T=\{e\})$ with a simply connected solvable group $F$
diffeomorphic to a Euclidean space and a torus $T$;
and $S$ is a Levi subgroup of $G$ centralized by $T$.
Note that we have $G= F \cdot K \,(F \cap K= \{e\})$, where $K=T\cdot S$ is
a reductive Lie group containing $H$.
\smallskip

In the case where $S$ is trivial, $G$ is a solvable Lie group; and an Inoue surface
of type $S^-$ is the only locally homogeneous l.c.K. manifold
with non-trivial $H$ (Example~1). In  the case where $S$ is non-trivial
we have ${1 \le \rm dim}\, R \le 2$ and $2 \le {\rm dim}\ S/H \le 3$; and thus
$R$ must be abelian. Hence we have $M=\R \times K/H$, where
$K/H$ is diffeomorphic to $U(2)/U(1)$ or $U(1,1)/U(1)$ for the case  ${\rm dim}\ R =2$,
which is a Hopf surface (Example 2) or a properly elliptic surface (\cite{W})
respectively;
and $K/H$ is diffeomorphic to $SU(2)$ or $\widetilde{SU}(1,1) \cong \widetilde{SL}(2,\R)$
for the case  ${\rm dim}\ R =1$ 
with $H=\{e\}$, which is to be excluded.
Note that the case $M=\R^2 \times K/H$ is excluded since it follows that
$K/H$ is diffeomorphic to $\C P^1$ or $\bf H$ (upper half plane); 
and $M$ could admit a K\"ahler structure, contradicting
that $\Gamma \backslash M$ admits no K\"ahler structures.
\end{proof}
\medskip

The classification of non-K\"ahler complex surfaces
except those of class VII with $b_2 > 0$ is known:
{\em Kodaira surfaces, 
Inoue surfaces, properly elliptic surfaces of odd type or Hopf surfaces}.
\smallskip

It is known (due to Wall \cite{W}) that a linear primary Hopf surface
with eigenvalues $\lambda_1, \lambda_2$ has a locally homogeneous
complex structure if and only if the condition $|\lambda_1| =|\lambda_2|$
holds. We will show in Example 2 that these Hopf surfaces admit locally
homogeneous l.c.K. structures. Therefore, we have obtained the following:
\smallskip

{\em Except for the class of diagonal Hopf surfaces
with eigenvalues $\lambda_1, \lambda_2 \,(|\lambda_1| \not=|\lambda_2|)$ (including
their secondary Hopf surfaces),
we can express, up to small deformations, each of l.c.K.
complex surfaces as $\Gamma \backslash G/H$,
where $G/H$ is a simply connected $4$-dimensional homogeneous l.c.K. manifold
and $\Gamma$ is a discrete subgroup of $G$}\,.

\section{Examples} 

\begin{Ex} 
{\rm An Inoue surface of type $S^-$ admits a locally homogeneous
l.c.K. structure of the form $\Gamma \backslash G/H$. In fact,
let $G=N \rtimes \R$ with the nilpotent Lie group $N$ defined by
$$ N=\{\left(
\begin{array}{ccc}
1 & x & z\\
0 & 1 & y\\
0 & 0 & 1
\end{array}
\right) | \; x \in \R,\, y, z \in \C
\},
$$
and the action $\phi(t): \R \rightarrow {\rm Aut}(N)$
defined by
$$\phi(t): 
\left(
\begin{array}{ccc}
1 & x & z\\
0 & 1 & y\\
0 & 0 & 1
\end{array}
\right)
\rightarrow
\left(
\begin{array}{ccc}
1 & a^t x & e^{\pi \im t}\, z\\
0 & 1 & a^{-t} e^{\pi \im t}\, y\\
0 & 0 & 1
\end{array}
\right),
$$
where $a, -\frac{1}{a} \,( a > 1)$ are real eigenvalues of
some $A \in GL(2, \Z)$.
If we take a closed subgroup $H$ of $N$ defined by
the condition $x=0, y=\im u, z=\im v \;(u, v \in \R)$ and
$\Gamma = \Lambda \rtimes \Z$ with $\Lambda$ a discrete subgroup of $N$
defined by the condition $x, y, z \in \Z$,
then $\phi(1)$ preserves $\Lambda$, and $\Gamma \backslash G/H$ is an Inoue
surface of type $S^-$ which has
an Inoue surface of type $S^+$ as a double covering
(\cite{Has} for more details). 
Note that the adjoint action of $H$ on $\g$ is trivial $\rm mod$ $\h$; hence
a l.c.K. structure on $\g/\h$ defines a homogeneous l.c.K. structure on $G/H$.
locally homogeneous manifold 
$\widehat{\Gamma} \backslash \widehat{G}$,
where $\widehat{G} =(N_\R \rtimes \Z_2) \rtimes \R$
with the action $\hat{\psi}: \Z_2 \rightarrow {\rm Aut} (N_\R)$ defined by
$$\hat{\psi}(s): 
\left(
\begin{array}{ccc}
1 & x & z\\
0 & 1 & y\\
0 & 0 & 1
\end{array}
\right)
\rightarrow
\left(
\begin{array}{ccc}
1 &  x &  (-1)^s z\\
0 & 1 &  (-1)^s y\\
0 & 0 & 1
\end{array}
\right),\; x,y,z \in \R,
$$
the action $\hat{\phi}: \R \rightarrow {\rm Aut}(N_\R)$ defined by
$$\hat{\phi}(t): 
\left(
\begin{array}{ccc}
1 & x & z\\
0 & 1 & y\\
0 & 0 & 1
\end{array}
\right)
\rightarrow
\left(
\begin{array}{ccc}
1 & a^t x &  z\\
0 & 1 & a^{-t} y\\
0 & 0 & 1
\end{array}
\right),\; x,y,z \in \R;
$$
and $\widehat{\Gamma} = (N_\Z \rtimes \Z_2) \rtimes \Z$
with the action $\tau: \Z \rightarrow {\rm Aut}(N_\Z \rtimes \Z_2)$
defined by $\tau(1)=\hat{\phi}(1) \hat{\psi}(1) \times 1$.
}
\end{Ex}

\begin{Ex} 
{\rm We can also consider $=S^1 \times S^3$ as a compact homogeneous
space $\tilde{G}/H$, where $\tilde{G}=S^1 \times U(2)$
with its Lie algebra 
$\tilde{\g}=\R \oplus {\mathfrak u}(2)$ and
$H=U(1)$ with its Lie algebra $\h$.
Then, we have a decomposition $\tilde{\g}=\m+\h$ for the subspace $\m$ of
$\tilde{\g}$ generated by $S, T, Y, Z$ and $\h$ generated by $W$, where
$$ 
S=\frac{1}{2}\left(
\begin{array}{cc}
\im & 0\\
0 & \im
\end{array}
\right),
\;
W=\frac{1}{2}\left(
\begin{array}{cc}
0 & 0\\
0 & \im
\end{array}
\right).
$$
Since we have $S=X+2 W$, we can take $\m'$ generated by
$T,X,Y,Z$ for $\m$; and homogeneous l.c.K. structures
on $\tilde{G}/H$ are the same as those on $G$. In other words
any homogeneous l.c.K. structures on $G$ can be extended as
those on $\tilde{G}/H$. 
\smallskip

Furthermore, we can construct locally homogeneous
l.c.K. manifolds $\Gamma \backslash \hat{G}/H$ for some discrete
subgroups $\Gamma$ of $\hat{G}$, where $\hat{G}=\R \times U(2)$.
For instance, let $\Gamma_{p,q}\, (p,q \not=0)$ be a discrete
subgroup of $\hat{G}$:
$$ 
\Gamma_{p,q}=\{(k,\left(
\begin{array}{cc}
e^{\im p k} & 0\\
0 & e^{\im q  k}
\end{array}
\right))
\in \R \times U(2)\, |\; k \in \Z
\}.
$$
Then $\Gamma_{p,q} \backslash \hat{G}/H$ is biholomorphic to
a Hopf surface $S_{p,q} =W/\Gamma_{\lambda_1,\lambda_2}$, where
$\Gamma_{\lambda_1,\lambda_2}$ is the cyclic group of automorphisms on
$W$ generated by

$$\phi: (z_1,z_2) \longrightarrow (\lambda_1 z_1, \lambda_2 z_2)$$
with $\lambda_1=e^{r+\im p}, \lambda_2=e^{r+\im q}, r \not=0.$
In fact, if we take a homogeneous complex structure $J_r$ on $\hat{G}/H$
induced from the diffeomorphism $\Phi_r: \hat{G}/H \rightarrow W$ defined by
$(t, z_1,z_2) \longrightarrow (e^{r t}  z_1, e^{r t} z_2)$, $\Phi_r$ induces
a biholomorphism between $\Gamma_{p,q} \backslash \hat{G}/H$ and $S_{p,q}$.
Note that in case $p=q$, $S_{p,q}$ is biholomorphic to $S_{\lambda}$ with
$\lambda=r+\im q$.
}
\end{Ex}
\medskip

We have an example of a compact locally homogeneous l.c.K. manifold of
non-compact reductive Lie group which is not of Vaisman type (\cite{ACHK}).

\begin{Ex} 
{\rm There exists a homogeneous l.c.K. structure on $\g=\R \oplus {\mathfrak sl}(2, \R)$,
which is not of Vaisman type. We can take a basis $\{W, X, Y, Z\}$ for $\g$ with
bracket multiplication defined by
$$[X,Y]=-Z,\, [Z,X]=Y,\, [Z,Y]=-X,$$
and all other brackets vanish.
We have a homogeneous complex
structure defined by
$$J Y=X, J X=-Y, J W=Z, J Z=-W,$$
and its compatible l.c.K. form $\Omega$ on $\g$ defined by
$$\Omega = z \wedge w + x \wedge y,$$
with the Lee form $\theta=w$, where $x, y, z, w$ are the Maurer-Cartan
forms corresponding to $X,Y,Z,W$ respectively. 
We can take another
l.c.K. form $\Omega_{\psi} = \psi \wedge w+ d \psi$, where 
$\psi=b y + c z \,(b, c \in \R)$ with $0 < b < c$ and $c^2-b^2=c$,
making the corresponding metric $h_{\psi}$ positive definite. The Lee field $\xi$ is
given as $\xi=\frac{1}{c^2-b^2} (cW+bX)$. It is easy to check that
$h([\xi, X],Y) + h(X, [\xi,Y]) \not\equiv 0;$
and thus $\xi$ is not a Killing field.
\smallskip

For any lattice $\Gamma$ of $G = \R\times {\widetilde{SL}(2, \R)}$ with
the above homogeneous l.c.K. structure,
we get a complex surface $\Gamma \backslash G$ (properly elliptic surface)
with locally homogeneous non-Vaisman l.c.K. structure.}
\end{Ex}

\bigskip

\address{Department of Mathematics\\
Faculty of Education\\
Niigata University\\
Niigata 950-2181\\
JAPAN}
{hasegawa@ed.niigata-u.ac.jp}
\address{
Department of Mathematics\\
Tokyo Metropolitan University\\
Hachioji 192-0397\\
JAPAN}
{kami@tmu.ac.jp}

\end{document}